\documentclass{amsart}
\usepackage{amsfonts}

\newtheorem{theorem}{Theorem}
\theoremstyle{plain}

\newtheorem{corollary}{Corollary}

\newtheorem{lemma}{Lemma}

\newtheorem{remark}{Remark}

\numberwithin{equation}{section}
\input{tcilatex}

\begin{document}
\title[Bombieri Inequality]{On the Bombieri Inequality in Inner Product Spaces}
\author{S.S. Dragomir}
\address{School of Computer Science and Mathematics\\
Victoria University of Technology\\
PO Box 14428, MCMC \\
Victoria 8001, Australia.}
\email{sever.dragomir@vu.edu.au}
\urladdr{http://rgmia.vu.edu.au/SSDragomirWeb.html}
\date{3 June, 2003.}
\subjclass{{26D15, 46C05.}}
\keywords{Bessel's inequality, Bombieri inequality.}

\begin{abstract}
New results related to the Bombieri generalisation of Bessel's inequality in
inner product spaces are given.
\end{abstract}

\maketitle

\section{Introduction}

Let $\left( H;\left( \cdot ,\cdot \right) \right) $ be an inner product
space over the real or complex number field $\mathbb{K}$. If $\left(
e_{i}\right) _{1\leq i\leq n}$ are orthonormal vectors in the inner product
space $H,$ i.e., $\left( e_{i},e_{j}\right) =\delta _{ij}$ for all $i,j\in
\left\{ 1,\dots ,n\right\} $ where $\delta _{ij}$ is the Kronecker delta,
then the following inequality is well known in the literature as Bessel's
inequality (see for example \cite[p. 391]{6b}): 
\begin{equation}
\sum_{i=1}^{n}\left| \left( x,e_{i}\right) \right| ^{2}\leq \left\|
x\right\| ^{2}\text{ \ for any }x\in H.  \label{1.1}
\end{equation}

For other results related to Bessel's inequality, see \cite{3b} -- \cite{5b}
and Chapter XV in the book \cite{6b}.

In 1971, E. Bombieri \cite{2ab} (see also \cite[p. 394]{6b}) gave the
following generalisation of Bessel's inequality.

\begin{theorem}
\label{t1.1}If $x,y_{1},\dots ,y_{n}$ are vectors in the inner product space 
$\left( H;\left( \cdot ,\cdot \right) \right) ,$ then the following
inequality: 
\begin{equation}
\sum_{i=1}^{n}\left| \left( x,y_{i}\right) \right| ^{2}\leq \left\|
x\right\| ^{2}\max_{1\leq i\leq n}\left\{ \sum_{j=1}^{n}\left| \left(
y_{i},y_{j}\right) \right| \right\} ,  \label{1.2}
\end{equation}
holds.
\end{theorem}

It is obvious that if $\left( y_{i}\right) _{1\leq i\leq n}$ are supposed to
be orthonormal, then from (\ref{1.2}) one would deduce Bessel's inequality (%
\ref{1.1}).

Another generalisation of Bessel's inequality was obtained by A. Selberg
(see for example \cite[p. 394]{6b}):

\begin{theorem}
\label{t1.2}Let $x,y_{1},\dots ,y_{n}$ be vectors in $H$ with $y_{i}\neq 0$ $%
\left( i=1,\dots ,n\right) .$ Then one has the inequality: 
\begin{equation}
\sum_{i=1}^{n}\frac{\left\vert \left( x,y_{i}\right) \right\vert ^{2}}{%
\sum_{j=1}^{n}\left\vert \left( y_{i},y_{j}\right) \right\vert }\leq
\left\Vert x\right\Vert ^{2}.  \label{1.3}
\end{equation}
\end{theorem}

In this case, also, if $\left( y_{i}\right) _{1\leq i\leq n}$ are
orthonormal, then from (\ref{1.3}) one may deduce Bessel's inequality.

Another type of inequality related to Bessel's result, was discovered in
1958 by H. Heilbronn \cite{5ab} (see also \cite[p. 395]{6b}).

\begin{theorem}
\label{t1.3}With the assumptions in Theorem \ref{t1.1}, one has 
\begin{equation}
\sum_{i=1}^{n}\left\vert \left( x,y_{i}\right) \right\vert \leq \left\Vert
x\right\Vert \left( \sum_{i,j=1}^{n}\left\vert \left( y_{i},y_{j}\right)
\right\vert \right) ^{\frac{1}{2}}.  \label{1.4}
\end{equation}
\end{theorem}

If in (\ref{1.4}) one chooses $y_{i}=e_{i}$ $\left( i=1,\dots ,n\right) ,$
where $\left( e_{i}\right) _{1\leq i\leq n}$ are orthonormal vectors in $H,$
then 
\begin{equation}
\sum_{i=1}^{n}\left| \left( x,e_{i}\right) \right| \leq \sqrt{n}\left\|
x\right\| ,\text{ \ for any \ }x\in H.  \label{1.5}
\end{equation}

In 1992 J.E. Pe\v{c}ari\'{c} \cite{7b} (see also \cite[p. 394]{6b}) proved
the following general inequality in inner product spaces.

\begin{theorem}
\label{t1.4}Let $x,y_{1},\dots ,y_{n}\in H$ and $c_{1},\dots ,c_{n}\in 
\mathbb{K}$. Then 
\begin{align}
\left| \sum_{i=1}^{n}c_{i}\left( x,y_{i}\right) \right| ^{2}& \leq \left\|
x\right\| ^{2}\sum_{i=1}^{n}\left| c_{i}\right| ^{2}\left(
\sum_{j=1}^{n}\left| \left( y_{i},y_{j}\right) \right| \right)  \label{1.6}
\\
& \leq \left\| x\right\| ^{2}\sum_{i=1}^{n}\left| c_{i}\right|
^{2}\max_{1\leq i\leq n}\left\{ \sum_{j=1}^{n}\left| \left(
y_{i},y_{j}\right) \right| \right\} .  \notag
\end{align}
\end{theorem}

He showed that the Bombieri inequality (\ref{1.2}) may be obtained from (\ref
{1.6}) for the choice $c_{i}=\overline{\left( x,y_{i}\right) }$ (using the
second inequality), the Selberg inequality (\ref{1.3}) may be obtained from
the first part of (\ref{1.6}) for the choice 
\begin{equation*}
c_{i}=\frac{\overline{\left( x,y_{i}\right) }}{\sum_{j=1}^{n}\left| \left(
y_{i},y_{j}\right) \right| },i\in \left\{ 1,\dots ,n\right\} ;
\end{equation*}
while the Heilbronn inequality (\ref{1.4}) may be obtained from the first
part of (\ref{1.6}) if one chooses $c_{i}=\frac{\overline{\left(
x,y_{i}\right) }}{\left| \left( x,y_{i}\right) \right| },$ for any $i\in
\left\{ 1,\dots ,n\right\} .$

For other results connected with the above ones, see \cite{4b} and \cite{5b}.

\section{Some Preliminary Results}

We start with the following lemma which is also interesting in itself.

\begin{lemma}
\label{l2.1}Let $z_{1},\dots ,z_{n}\in H$ and $\alpha _{1},\dots ,\alpha
_{n}\in \mathbb{K}.$ Then one has the inequality: 
\begin{equation}
\left\| \sum_{i=1}^{n}\alpha _{i}z_{i}\right\| ^{2}  \label{2.1}
\end{equation}
\begin{equation*}
\leq \left\{ 
\begin{array}{l}
\max\limits_{1\leq k\leq n}\left| \alpha _{k}\right|
^{2}\sum\limits_{i,j=1}^{n}\left| \left( z_{i},z_{j}\right) \right| ; \\ 
\\ 
\max\limits_{1\leq k\leq n}\left| \alpha _{k}\right| \left(
\sum\limits_{i=1}^{n}\left| \alpha _{i}\right| ^{r}\right) ^{\frac{1}{r}%
}\left( \sum\limits_{i=1}^{n}\left( \sum\limits_{j=1}^{n}\left| \left(
z_{i},z_{j}\right) \right| \right) ^{s}\right) ^{\frac{1}{s}},\ \ \ r>1,\ 
\frac{1}{r}+\frac{1}{s}=1; \\ 
\\ 
\max\limits_{1\leq k\leq n}\left| \alpha _{k}\right|
\sum\limits_{k=1}^{n}\left| \alpha _{k}\right| \max\limits_{1\leq i\leq
n}\left( \sum\limits_{j=1}^{n}\left| \left( z_{i},z_{j}\right) \right|
\right) ; \\ 
\\ 
\left( \sum\limits_{k=1}^{n}\left| \alpha _{k}\right| ^{p}\right) ^{\frac{1}{%
p}}\max\limits_{1\leq i\leq n}\left| \alpha _{i}\right| \left(
\sum\limits_{i=1}^{n}\left( \sum\limits_{j=1}^{n}\left| \left(
z_{i},z_{j}\right) \right| \right) ^{q}\right) ^{\frac{1}{q}},\ \ \ p>1,\ 
\frac{1}{p}+\frac{1}{q}=1; \\ 
\\ 
\left( \sum\limits_{k=1}^{n}\left| \alpha _{k}\right| ^{p}\right) ^{\frac{1}{%
p}}\left( \sum\limits_{i=1}^{n}\left| \alpha _{i}\right| ^{t}\right) ^{\frac{%
1}{t}}\left[ \sum\limits_{i=1}^{n}\left( \sum\limits_{j=1}^{n}\left| \left(
z_{i},z_{j}\right) \right| ^{q}\right) ^{\frac{u}{q}}\right] ^{\frac{1}{u}%
},\ \ \ p>1,\ \frac{1}{p}+\frac{1}{q}=1; \\ 
\hfill t>1,\ \frac{1}{t}+\frac{1}{u}=1; \\ 
\\ 
\left( \sum\limits_{k=1}^{n}\left| \alpha _{k}\right| ^{p}\right) ^{\frac{1}{%
p}}\sum\limits_{i=1}^{n}\left| \alpha _{i}\right| \max\limits_{1\leq i\leq
n}\left\{ \left( \sum\limits_{j=1}^{n}\left| \left( z_{i},z_{j}\right)
\right| ^{q}\right) ^{\frac{1}{q}}\right\} ,\ \ \ p>1,\ \frac{1}{p}+\frac{1}{%
q}=1; \\ 
\\ 
\sum\limits_{k=1}^{n}\left| \alpha _{k}\right| \ \max\limits_{1\leq i\leq
n}\left| \alpha _{i}\right| \ \sum\limits_{i=1}^{n}\left[ \max\limits_{1\leq
j\leq n}\left| \left( z_{i},z_{j}\right) \right| \right] ; \\ 
\\ 
\sum\limits_{k=1}^{n}\left| \alpha _{k}\right| \ \left(
\sum\limits_{i=1}^{n}\left| \alpha _{i}\right| ^{m}\right) ^{\frac{1}{m}%
}\left( \sum\limits_{i=1}^{n}\left[ \max\limits_{1\leq j\leq n}\left| \left(
z_{i},z_{j}\right) \right| \right] ^{l}\right) ^{\frac{1}{l}},\ \ \ m>1,\ 
\frac{1}{m}+\frac{1}{l}=1; \\ 
\\ 
\left( \sum\limits_{k=1}^{n}\left| \alpha _{k}\right| \right)
^{2}\max\limits_{i,1\leq j\leq n}\left| \left( z_{i},z_{j}\right) \right| .
\end{array}
\right.
\end{equation*}
\end{lemma}

\begin{proof}
We observe that 
\begin{align}
\left\| \sum_{i=1}^{n}\alpha _{i}z_{i}\right\| ^{2}& =\left(
\sum_{i=1}^{n}\alpha _{i}z_{i},\sum_{j=1}^{n}\alpha _{j}z_{j}\right)
\label{2.2} \\
& =\sum_{i=1}^{n}\sum_{j=1}^{n}\alpha _{i}\overline{\alpha _{j}}\left(
z_{i},z_{j}\right) =\left| \sum_{i=1}^{n}\sum_{j=1}^{n}\alpha _{i}\overline{%
\alpha _{j}}\left( z_{i},z_{j}\right) \right|  \notag \\
& \leq \sum_{i=1}^{n}\sum_{j=1}^{n}\left| \alpha _{i}\right| \left| \alpha
_{j}\right| \left| \left( z_{i},z_{j}\right) \right| =\sum_{i=1}^{n}\left|
\alpha _{i}\right| \left( \sum_{j=1}^{n}\left| \alpha _{j}\right| \left|
\left( z_{i},z_{j}\right) \right| \right)  \notag \\
& :=M.  \notag
\end{align}
Using H\"{o}lder's inequality, we may write that 
\begin{equation}
\sum_{j=1}^{n}\left| \alpha _{j}\right| \left| \left( z_{i},z_{j}\right)
\right| \leq \left\{ 
\begin{array}{l}
\max\limits_{1\leq k\leq n}\left| \alpha _{k}\right|
\sum\limits_{j=1}^{n}\left| \left( z_{i},z_{j}\right) \right| \\ 
\\ 
\left( \sum\limits_{k=1}^{n}\left| \alpha _{k}\right| ^{p}\right) ^{\frac{1}{%
p}}\left( \sum\limits_{j=1}^{n}\left| \left( z_{i},z_{j}\right) \right|
^{q}\right) ^{\frac{1}{q}},\ \ p>1,\ \frac{1}{p}+\frac{1}{q}=1; \\ 
\\ 
\sum\limits_{k=1}^{n}\left| \alpha _{k}\right| \ \max\limits_{1\leq j\leq
n}\left| \left( z_{i},z_{j}\right) \right|
\end{array}
\right.  \label{2.3}
\end{equation}
for any $i\in \left\{ 1,\dots ,n\right\} ,$ giving 
\begin{equation}
M\leq \left\{ 
\begin{array}{l}
\max\limits_{1\leq k\leq n}\left| \alpha _{k}\right|
\sum\limits_{i=1}^{n}\left| \alpha _{i}\right| \sum\limits_{j=1}^{n}\left|
\left( z_{i},z_{j}\right) \right| =:M_{1}; \\ 
\\ 
\left( \sum\limits_{k=1}^{n}\left| \alpha _{k}\right| ^{p}\right) ^{\frac{1}{%
p}}\sum\limits_{i=1}^{n}\left| \alpha _{i}\right| \left(
\sum\limits_{j=1}^{n}\left| \left( z_{i},z_{j}\right) \right| ^{q}\right) ^{%
\frac{1}{q}}:=M_{p},\  \\ 
\hfill p>1,\ \frac{1}{p}+\frac{1}{q}=1; \\ 
\sum\limits_{k=1}^{n}\left| \alpha _{k}\right| \sum\limits_{i=1}^{n}\left|
\alpha _{i}\right| \ \max\limits_{1\leq j\leq n}\left| \left(
z_{i},z_{j}\right) \right| =:M_{\infty }.
\end{array}
\right.  \label{2.4}
\end{equation}
By H\"{o}lder's inequality we also have: 
\begin{multline}
\sum\limits_{i=1}^{n}\left| \alpha _{i}\right| \left(
\sum\limits_{j=1}^{n}\left| \left( z_{i},z_{j}\right) \right| \right)
\label{2.5} \\
\leq \left\{ 
\begin{array}{l}
\max\limits_{1\leq i\leq n}\left| \alpha _{i}\right|
\sum\limits_{i,j=1}^{n}\left| \left( z_{i},z_{j}\right) \right| ; \\ 
\\ 
\left( \sum\limits_{i=1}^{n}\left| \alpha _{i}\right| ^{r}\right) ^{\frac{1}{%
r}}\left( \sum\limits_{i=1}^{n}\left( \sum\limits_{j=1}^{n}\left| \left(
z_{i},z_{j}\right) \right| \right) ^{s}\right) ^{\frac{1}{s}},\ \ \ r>1,\ 
\frac{1}{r}+\frac{1}{s}=1; \\ 
\\ 
\sum\limits_{i=1}^{n}\left| \alpha _{i}\right| \ \max\limits_{1\leq i\leq
n}\left( \sum\limits_{j=1}^{n}\left| \left( z_{i},z_{j}\right) \right|
\right) ;
\end{array}
\right.
\end{multline}
and thus 
\begin{equation*}
M_{1}\leq \left\{ 
\begin{array}{l}
\max\limits_{1\leq k\leq n}\left| \alpha _{k}\right|
^{2}\sum\limits_{i,j=1}^{n}\left| \left( z_{i},z_{j}\right) \right| ; \\ 
\\ 
\max\limits_{1\leq k\leq n}\left| \alpha _{k}\right| \left(
\sum\limits_{i=1}^{n}\left| \alpha _{i}\right| ^{r}\right) ^{\frac{1}{r}%
}\left( \sum\limits_{i=1}^{n}\left( \sum\limits_{j=1}^{n}\left| \left(
z_{i},z_{j}\right) \right| \right) ^{s}\right) ^{\frac{1}{s}},\ \ \ r>1,\ 
\frac{1}{r}+\frac{1}{s}=1; \\ 
\\ 
\max\limits_{1\leq k\leq n}\left| \alpha _{k}\right|
\sum\limits_{i=1}^{n}\left| \alpha _{i}\right| \ \max\limits_{1\leq i\leq
n}\left( \sum\limits_{j=1}^{n}\left| \left( z_{i},z_{j}\right) \right|
\right) ;
\end{array}
\right.
\end{equation*}
and the first 3 inequalities in (\ref{2.1}) are obtained.

By H\"{o}lder's inequality we also have: 
\begin{multline*}
M_{p}\leq \left( \sum\limits_{k=1}^{n}\left| \alpha _{k}\right| ^{p}\right)
^{\frac{1}{p}} \\
\times \left\{ 
\begin{array}{l}
\max\limits_{1\leq i\leq n}\left| \alpha _{i}\right|
\sum\limits_{i=1}^{n}\left( \sum\limits_{j=1}^{n}\left| \left(
z_{i},z_{j}\right) \right| ^{q}\right) ^{\frac{1}{q}}; \\ 
\\ 
\left( \sum\limits_{i=1}^{n}\left| \alpha _{i}\right| ^{t}\right) ^{\frac{1}{%
t}}\left( \sum\limits_{i=1}^{n}\left( \sum\limits_{j=1}^{n}\left| \left(
z_{i},z_{j}\right) \right| ^{q}\right) ^{\frac{u}{q}}\right) ^{\frac{1}{u}%
},\ \ \ \hfill t>1,\ \frac{1}{t}+\frac{1}{u}=1; \\ 
\\ 
\sum\limits_{i=1}^{n}\left| \alpha _{i}\right| \ \max\limits_{1\leq i\leq
n}\left\{ \left( \sum\limits_{j=1}^{n}\left| \left( z_{i},z_{j}\right)
\right| ^{q}\right) ^{\frac{1}{q}}\right\} ;
\end{array}
\right.
\end{multline*}
and the next 3 inequalities in (\ref{2.1}) are proved.

Finally, by the same H\"{o}lder inequality we may state that: 
\begin{equation*}
M_{\infty }\leq \sum\limits_{k=1}^{n}\left| \alpha _{k}\right| \times
\left\{ 
\begin{array}{l}
\max\limits_{1\leq i\leq n}\left| \alpha _{i}\right|
\sum\limits_{i=1}^{n}\left( \max\limits_{1\leq j\leq n}\left| \left(
z_{i},z_{j}\right) \right| \right) ; \\ 
\\ 
\left( \sum\limits_{i=1}^{n}\left| \alpha _{i}\right| ^{m}\right) ^{\frac{1}{%
m}}\left( \sum\limits_{i=1}^{n}\left( \max\limits_{1\leq j\leq n}\left|
\left( z_{i},z_{j}\right) \right| \right) ^{l}\right) ^{\frac{1}{l}},\ \ \
\hfill m>1,\ \frac{1}{m}+\frac{1}{l}=1; \\ 
\\ 
\sum\limits_{i=1}^{n}\left| \alpha _{i}\right| \ \max\limits_{1\leq i,j\leq
n}\left| \left( z_{i},z_{j}\right) \right| ;
\end{array}
\right.
\end{equation*}
and the last 3 inequalities in (\ref{2.1}) are proved.
\end{proof}

If we would like to have some bounds for $\left\| \sum_{i=1}^{n}\alpha
_{i}z_{i}\right\| ^{2}$ in terms of $\sum_{i=1}^{n}\left| \alpha _{i}\right|
^{2},$ then the following corollaries may be used.

\begin{corollary}
\label{c2.1.a}Let $z_{1},\dots ,z_{n}$ and $\alpha _{1},\dots ,\alpha _{n}$
be as in Lemma \ref{l2.1}. If $1<p\leq 2$, $1<t\leq 2,$ then one has the
inequality 
\begin{equation}
\left\Vert \sum_{i=1}^{n}\alpha _{i}z_{i}\right\Vert ^{2}\leq n^{\frac{1}{p}+%
\frac{1}{t}-1}\sum\limits_{k=1}^{n}\left\vert \alpha _{k}\right\vert ^{2}%
\left[ \sum\limits_{i=1}^{n}\left( \sum\limits_{j=1}^{n}\left\vert \left(
z_{i},z_{j}\right) \right\vert ^{q}\right) ^{\frac{u}{q}}\right] ^{\frac{1}{u%
}}  \label{2.5.a}
\end{equation}
where $\frac{1}{p}+\frac{1}{q}=1,$ $\frac{1}{t}+\frac{1}{u}=1.$
\end{corollary}

\begin{proof}
Observe, by the monotonicity of power means, we may write that 
\begin{align*}
\left( \frac{\sum_{k=1}^{n}\left\vert \alpha _{k}\right\vert ^{p}}{n}\right)
^{\frac{1}{p}}& \leq \left( \frac{\sum_{k=1}^{n}\left\vert \alpha
_{k}\right\vert ^{2}}{n}\right) ^{\frac{1}{2}};\ \ 1<p\leq 2, \\
\left( \frac{\sum_{k=1}^{n}\left\vert \alpha _{k}\right\vert ^{t}}{n}\right)
^{\frac{1}{t}}& \leq \left( \frac{\sum_{k=1}^{n}\left\vert \alpha
_{k}\right\vert ^{2}}{n}\right) ^{\frac{1}{2}};\ \ 1<t\leq 2,
\end{align*}
from where we get 
\begin{align*}
\left( \sum_{k=1}^{n}\left\vert \alpha _{k}\right\vert ^{p}\right) ^{\frac{1%
}{p}}& \leq n^{\frac{1}{p}-\frac{1}{2}}\left( \sum_{k=1}^{n}\left\vert
\alpha _{k}\right\vert ^{2}\right) ^{\frac{1}{2}}, \\
\left( \sum_{k=1}^{n}\left\vert \alpha _{k}\right\vert ^{t}\right) ^{\frac{1%
}{t}}& \leq n^{\frac{1}{t}-\frac{1}{2}}\left( \sum_{k=1}^{n}\left\vert
\alpha _{k}\right\vert ^{2}\right) ^{\frac{1}{2}}.
\end{align*}
Using the fifth inequality in (\ref{2.1}), we then deduce (\ref{2.5.a}).
\end{proof}

\begin{remark}
\label{r2.1.b}An interesting particular case is the one for $p=q=t=u=2,$
giving 
\begin{equation}
\left\Vert \sum_{i=1}^{n}\alpha _{i}z_{i}\right\Vert ^{2}\leq
\sum_{k=1}^{n}\left\vert \alpha _{k}\right\vert ^{2}\left(
\sum\limits_{i,j=1}^{n}\left\vert \left( z_{i},z_{j}\right) \right\vert
^{2}\right) ^{\frac{1}{2}}.  \label{2.5.b}
\end{equation}
\end{remark}

\begin{corollary}
\label{c2.1.c}With the assumptions of Lemma \ref{l2.1} and if $1<p\leq 2,$
then 
\begin{equation}
\left\| \sum_{i=1}^{n}\alpha _{i}z_{i}\right\| ^{2}\leq n^{\frac{1}{p}%
}\sum_{k=1}^{n}\left| \alpha _{k}\right| ^{2}\max\limits_{1\leq i\leq n}%
\left[ \left( \sum\limits_{j=1}^{n}\left| \left( z_{i},z_{j}\right) \right|
^{q}\right) ^{\frac{1}{q}}\right] ,  \label{2.5.c}
\end{equation}
where $\frac{1}{p}+\frac{1}{q}=1.$
\end{corollary}

\begin{proof}
Since 
\begin{equation*}
\left( \sum_{k=1}^{n}\left\vert \alpha _{k}\right\vert ^{p}\right) ^{\frac{1%
}{p}}\leq n^{\frac{1}{p}-\frac{1}{2}}\left( \sum_{k=1}^{n}\left\vert \alpha
_{k}\right\vert ^{2}\right) ^{\frac{1}{2}},
\end{equation*}
and 
\begin{equation*}
\sum_{k=1}^{n}\left\vert \alpha _{k}\right\vert \leq n^{\frac{1}{2}}\left(
\sum_{k=1}^{n}\left\vert \alpha _{k}\right\vert ^{2}\right) ^{\frac{1}{2}},
\end{equation*}
then by the sixth inequality in (\ref{2.1}) we deduce (\ref{2.5.c}).
\end{proof}

In a similar fashion, one may prove the following two corollaries.

\begin{corollary}
\label{c2.1.1}With the assumptions of Lemma \ref{l2.1} and if $1<m\leq 2,$
then 
\begin{equation}
\left\| \sum_{i=1}^{n}\alpha _{i}z_{i}\right\| ^{2}\leq n^{\frac{1}{m}%
}\sum_{k=1}^{n}\left| \alpha _{k}\right| ^{2}\left( \sum\limits_{i=1}^{n}%
\left[ \max\limits_{1\leq j\leq n}\left| \left( z_{i},z_{j}\right) \right| %
\right] ^{l}\right) ^{\frac{1}{l}},  \label{2.5.d}
\end{equation}
where $\frac{1}{m}+\frac{1}{l}=1.$
\end{corollary}

\begin{corollary}
\label{c2.1.e}With the assumptions of Lemma \ref{l2.1}, we have: 
\begin{equation}
\left\| \sum_{i=1}^{n}\alpha _{i}z_{i}\right\| ^{2}\leq
n\sum_{k=1}^{n}\left| \alpha _{k}\right| ^{2}\ \max\limits_{1\leq i,j\leq
n}\left| \left( z_{i},z_{j}\right) \right| .  \label{2.5.e}
\end{equation}
\end{corollary}

The following lemma may be of interest as well.

\begin{lemma}
\label{l2.2}With the assumptions of Lemma \ref{l2.1}, one has the
inequalities 
\begin{equation}
\left\| \sum_{i=1}^{n}\alpha _{i}z_{i}\right\| ^{2}\leq \sum_{i=1}^{n}\left|
\alpha _{i}\right| ^{2}\sum\limits_{j=1}^{n}\left| \left( z_{i},z_{j}\right)
\right|  \label{2.6}
\end{equation}
\begin{equation*}
\leq \left\{ 
\begin{array}{l}
\sum\limits_{i=1}^{n}\left| \alpha _{i}\right| ^{2}\ \max\limits_{1\leq
i\leq n}\left[ \sum\limits_{j=1}^{n}\left| \left( z_{i},z_{j}\right) \right| 
\right] ; \\ 
\\ 
\left( \sum\limits_{i=1}^{n}\left| \alpha _{i}\right| ^{2p}\right) ^{\frac{1%
}{p}}\left( \left( \sum\limits_{j=1}^{n}\left| \left( z_{i},z_{j}\right)
\right| \right) ^{q}\right) ^{\frac{1}{q}},\ \ \ \ \ \ p>1,\ \frac{1}{p}+%
\frac{1}{q}=1; \\ 
\\ 
\max\limits_{1\leq i\leq n}\left| \alpha _{i}\right| ^{2}\
\sum\limits_{i,j=1}^{n}\left| \left( z_{i},z_{j}\right) \right| .
\end{array}
\right.
\end{equation*}
\end{lemma}

\begin{proof}
As in Lemma \ref{l2.1}, we know that 
\begin{equation}
\left\Vert \sum_{i=1}^{n}\alpha _{i}z_{i}\right\Vert ^{2}\leq
\sum_{i=1}^{n}\sum_{j=1}^{n}\left\vert \alpha _{i}\right\vert \left\vert
\alpha _{j}\right\vert \left\vert \left( z_{i},z_{j}\right) \right\vert .
\label{2.7}
\end{equation}
Using the simple observation that (see also \cite[p. 394]{6b}) 
\begin{equation*}
\left\vert \alpha _{i}\right\vert \left\vert \alpha _{j}\right\vert \leq 
\frac{1}{2}\left( \left\vert \alpha _{i}\right\vert ^{2}+\left\vert \alpha
_{j}\right\vert ^{2}\right) ,\ \ \ i,j\in \left\{ 1,\dots ,n\right\}
\end{equation*}
we have 
\begin{align*}
\sum_{i=1}^{n}\sum_{j=1}^{n}\left\vert \alpha _{i}\right\vert \left\vert
\alpha _{j}\right\vert \left\vert \left( z_{i},z_{j}\right) \right\vert &
\leq \frac{1}{2}\sum\limits_{i,j=1}^{n}\left( \left\vert \alpha
_{i}\right\vert ^{2}+\left\vert \alpha _{j}\right\vert ^{2}\right)
\left\vert \left( z_{i},z_{j}\right) \right\vert \\
& =\frac{1}{2}\left[ \sum\limits_{i,j=1}^{n}\left\vert \alpha
_{i}\right\vert ^{2}\left\vert \left( z_{i},z_{j}\right) \right\vert
+\sum\limits_{i,j=1}^{n}\left\vert \alpha _{j}\right\vert ^{2}\left\vert
\left( z_{i},z_{j}\right) \right\vert \right] \\
& =\sum\limits_{i,j=1}^{n}\left\vert \alpha _{i}\right\vert ^{2}\left\vert
\left( z_{i},z_{j}\right) \right\vert ,
\end{align*}
which proves the first inequality in (\ref{2.6}).

The second part follows by H\"{o}lder's inequality and we omit the details.
\end{proof}

\begin{remark}
The first part in (\ref{2.6}) is the inequality obtained by Pe\v{c}ari\'{c}
in \cite{7b}.
\end{remark}

\section{Some Pe\v{c}ari\'{c} Type Inequalities}

We are now able to point out the following result which complements the
inequality (\ref{1.6}) due to J.E. Pe\v{c}ari\'{c} \cite{7b} (see also 
\cite[p. 394]{6b}).

\begin{theorem}
\label{t3.1}Let $x,y_{1},\dots ,y_{n}$ be vectors of an inner product space $%
\left( H;\left( \cdot ,\cdot \right) \right) $ and $c_{1},\dots ,c_{n}\in 
\mathbb{K}$. Then one has the inequalities: 
\begin{equation}
\left| \sum\limits_{i=1}^{n}c_{i}\left( x,y_{i}\right) \right| ^{2}
\label{3.1}
\end{equation}
\begin{equation*}
\leq \left\| x\right\| ^{2}\times \left\{ 
\begin{array}{l}
\max\limits_{1\leq k\leq n}\left| c_{k}\right|
^{2}\sum\limits_{i,j=1}^{n}\left| \left( y_{i},y_{j}\right) \right| ; \\ 
\\ 
\max\limits_{1\leq k\leq n}\left| c_{k}\right| \left(
\sum\limits_{i=1}^{n}\left| c_{i}\right| ^{r}\right) ^{\frac{1}{r}}\left[
\sum\limits_{i=1}^{n}\left( \sum\limits_{j=1}^{n}\left| \left(
y_{i},y_{j}\right) \right| \right) ^{s}\right] ^{\frac{1}{s}},\ \ \ r>1,\ 
\frac{1}{r}+\frac{1}{s}=1; \\ 
\\ 
\max\limits_{1\leq k\leq n}\left| c_{k}\right| \sum\limits_{k=1}^{n}\left|
c_{k}\right| \max\limits_{1\leq i\leq n}\left( \sum\limits_{j=1}^{n}\left|
\left( y_{i},y_{j}\right) \right| \right) ; \\ 
\\ 
\left( \sum\limits_{k=1}^{n}\left| c_{k}\right| ^{p}\right) ^{\frac{1}{p}%
}\max\limits_{1\leq i\leq n}\left| c_{i}\right| \left(
\sum\limits_{i=1}^{n}\left( \sum\limits_{j=1}^{n}\left| \left(
y_{i},y_{j}\right) \right| \right) ^{q}\right) ^{\frac{1}{q}},\ \ \ p>1,\ 
\frac{1}{p}+\frac{1}{q}=1; \\ 
\\ 
\left( \sum\limits_{k=1}^{n}\left| c_{k}\right| ^{p}\right) ^{\frac{1}{p}%
}\left( \sum\limits_{i=1}^{n}\left| c_{i}\right| ^{t}\right) ^{\frac{1}{t}}%
\left[ \sum\limits_{i=1}^{n}\left( \sum\limits_{j=1}^{n}\left| \left(
y_{i},y_{j}\right) \right| ^{q}\right) ^{\frac{u}{q}}\right] ^{\frac{1}{u}%
},\ \ \ p>1,\ \frac{1}{p}+\frac{1}{q}=1; \\ 
\hfill t>1,\ \frac{1}{t}+\frac{1}{u}=1; \\ 
\\ 
\left( \sum\limits_{k=1}^{n}\left| c_{k}\right| ^{p}\right) ^{\frac{1}{p}%
}\sum\limits_{i=1}^{n}\left| c_{i}\right| \max\limits_{1\leq i\leq n}\left\{
\left( \sum\limits_{j=1}^{n}\left| \left( y_{i},y_{j}\right) \right|
^{q}\right) ^{\frac{1}{q}}\right\} ,\ \ \ p>1,\ \frac{1}{p}+\frac{1}{q}=1;
\\ 
\\ 
\sum\limits_{k=1}^{n}\left| c_{k}\right| \ \max\limits_{1\leq i\leq n}\left|
c_{i}\right| \ \sum\limits_{i=1}^{n}\left[ \max\limits_{1\leq j\leq n}\left|
\left( y_{i},y_{j}\right) \right| \right] ; \\ 
\\ 
\sum\limits_{k=1}^{n}\left| c_{k}\right| \ \left(
\sum\limits_{i=1}^{n}\left| c_{i}\right| ^{m}\right) ^{\frac{1}{m}}\left(
\sum\limits_{i=1}^{n}\left[ \max\limits_{1\leq j\leq n}\left| \left(
y_{i},y_{j}\right) \right| \right] ^{l}\right) ^{\frac{1}{l}},\ \ \ m>1,\ 
\frac{1}{m}+\frac{1}{l}=1; \\ 
\\ 
\left( \sum\limits_{k=1}^{n}\left| c_{k}\right| \right)
^{2}\max\limits_{i,1\leq j\leq n}\left| \left( y_{i},y_{j}\right) \right| .
\end{array}
\right.
\end{equation*}
\end{theorem}

\begin{proof}
We note that 
\begin{equation*}
\sum\limits_{i=1}^{n}c_{i}\left( x,y_{i}\right) =\left(
x,\sum\limits_{i=1}^{n}\overline{c_{i}}y_{i}\right) .
\end{equation*}
Using Schwarz's inequality in inner product spaces, we have 
\begin{equation}
\left| \sum\limits_{i=1}^{n}c_{i}\left( x,y_{i}\right) \right| ^{2}\leq
\left\| x\right\| ^{2}\left\| \sum\limits_{i=1}^{n}\overline{c_{i}}%
y_{i}\right\| ^{2}.  \label{3.2}
\end{equation}
Finally, using Lemma \ref{l2.1} with $\alpha _{i}=\overline{c_{i}},$ $%
z_{i}=y_{i}$ $\left( i=1,\dots ,n\right) ,$ we deduce the desired inequality
(\ref{3.1}). We omit the details.
\end{proof}

The following corollaries may be useful if one needs bounds in terms of $%
\sum_{i=1}^{n}\left| c_{i}\right| ^{2}.$

\begin{corollary}
\label{c3.2}With the assumptions in Theorem \ref{t3.1} and if $1<p\leq 2,$ $%
1<t\leq 2,$ $\frac{1}{p}+\frac{1}{q}=1,$ $\frac{1}{t}+\frac{1}{u}=1,$ one
has the inequality: 
\begin{equation}
\left\vert \sum\limits_{i=1}^{n}c_{i}\left( x,y_{i}\right) \right\vert
^{2}\leq \left\Vert x\right\Vert ^{2}n^{\frac{1}{p}+\frac{1}{t}%
-1}\sum\limits_{i=1}^{n}\left\vert c_{i}\right\vert ^{2}\left[
\sum\limits_{i=1}^{n}\left( \sum\limits_{j=1}^{n}\left\vert \left(
y_{i},y_{j}\right) \right\vert ^{q}\right) ^{\frac{u}{q}}\right] ^{\frac{1}{u%
}},  \label{3.3}
\end{equation}
and, in particular, for $p=q=t=u=2,$%
\begin{equation}
\left\vert \sum\limits_{i=1}^{n}c_{i}\left( x,y_{i}\right) \right\vert
^{2}\leq \left\Vert x\right\Vert ^{2}\sum\limits_{i=1}^{n}\left\vert
c_{i}\right\vert ^{2}\left( \sum\limits_{i,j=1}^{n}\left\vert \left(
y_{i},y_{j}\right) \right\vert ^{2}\right) ^{\frac{1}{2}}.  \label{3.4}
\end{equation}
\end{corollary}

The proof is similar to the one in Corollary \ref{c2.1.a}.

\begin{corollary}
\label{c3.3}With the assumptions in Theorem \ref{t3.1} and if $1<p\leq 2,$
then 
\begin{equation}
\left| \sum\limits_{i=1}^{n}c_{i}\left( x,y_{i}\right) \right| ^{2}\leq
\left\| x\right\| ^{2}n^{\frac{1}{p}}\sum\limits_{k=1}^{n}\left|
c_{k}\right| ^{2}\ \max\limits_{1\leq i\leq n}\left[ \sum\limits_{j=1}^{n}%
\left| \left( y_{i},y_{j}\right) \right| ^{q}\right] ^{\frac{1}{q}},
\label{3.5}
\end{equation}
where $\frac{1}{p}+\frac{1}{q}=1.$
\end{corollary}

The proof is similar to the one in Corollary \ref{c2.1.c}.

The following two inequalities also hold.

\begin{corollary}
\label{c3.4}With the above assumptions for $X,y_{i},c_{i}$ and if $1<m\leq
2, $ then 
\begin{equation}
\left| \sum\limits_{i=1}^{n}c_{i}\left( x,y_{i}\right) \right| ^{2}\leq
\left\| x\right\| ^{2}n^{\frac{1}{m}}\sum\limits_{k=1}^{n}\left|
c_{k}\right| ^{2}\left( \sum\limits_{i=1}^{n}\left[ \max\limits_{1\leq j\leq
n}\left| \left( y_{i},y_{j}\right) \right| \right] ^{l}\right) ^{\frac{1}{l}%
},  \label{3.6}
\end{equation}
where $\frac{1}{m}+\frac{1}{l}=1.$
\end{corollary}

\begin{corollary}
\label{c3.5}With the above assumptions for $X,y_{i},c_{i},$ one has 
\begin{equation}
\left| \sum\limits_{i=1}^{n}c_{i}\left( x,y_{i}\right) \right| ^{2}\leq
\left\| x\right\| ^{2}n\sum\limits_{k=1}^{n}\left| c_{k}\right|
^{2}\max\limits_{i,1\leq j\leq n}\left| \left( y_{i},y_{j}\right) \right| .
\label{3.7}
\end{equation}
\end{corollary}

Using Lemma \ref{l2.2}, we may state the following result as well.

\begin{remark}
\label{t3.6}With the assumptions of Theorem \ref{t3.1}, one has the
inequalities: 
\begin{equation}
\left| \sum\limits_{i=1}^{n}c_{i}\left( x,y_{i}\right) \right| ^{2}\leq
\left\| x\right\| ^{2}\sum\limits_{i=1}^{n}\left| c_{i}\right|
^{2}\sum\limits_{j=1}^{n}\left| \left( y_{i},y_{j}\right) \right|
\label{3.8}
\end{equation}
\begin{equation*}
\leq \left\| x\right\| ^{2}\times \left\{ 
\begin{array}{l}
\sum\limits_{i=1}^{n}\left| c_{i}\right| ^{2}\ \max\limits_{1\leq i\leq n} 
\left[ \sum\limits_{j=1}^{n}\left| \left( y_{i},y_{j}\right) \right| \right]
; \\ 
\\ 
\left( \sum\limits_{i=1}^{n}\left| c_{i}\right| ^{2p}\right) ^{\frac{1}{p}%
}\left( \sum\limits_{i=1}^{n}\left( \sum\limits_{j=1}^{n}\left| \left(
y_{i},y_{j}\right) \right| \right) ^{q}\right) ^{\frac{1}{q}},\ \ p>1,\ 
\frac{1}{p}+\frac{1}{q}=1; \\ 
\\ 
\max\limits_{1\leq i\leq n}\left| c_{i}\right| ^{2}\
\sum\limits_{i,j=1}^{n}\left| \left( y_{i},y_{j}\right) \right| ;
\end{array}
\right.
\end{equation*}
that provide some alternatives to Pe\v{c}ari\'{c}'s result (\ref{1.6}).
\end{remark}

\section{Some Inequalities of Bombieri Type}

In this section we point out some inequalities of Bombieri type that may be
obtained from (\ref{3.1}) on choosing $c_{i}=\overline{\left( x,y_{i}\right) 
}$ $\left( i=1,\dots ,n\right) .$

If the above choice was made in the first inequality in (\ref{3.1}), then
one would obtain: 
\begin{equation*}
\left( \sum\limits_{i=1}^{n}\left| \left( x,y_{i}\right) \right| ^{2}\right)
^{2}\leq \left\| x\right\| ^{2}\max\limits_{1\leq i\leq n}\left| \left(
x,y_{i}\right) \right| ^{2}\sum\limits_{i,j=1}^{n}\left| \left(
y_{i},y_{j}\right) \right|
\end{equation*}
giving, by taking the square root, 
\begin{equation}
\sum\limits_{i=1}^{n}\left| \left( x,y_{i}\right) \right| ^{2}\leq \left\|
x\right\| \max\limits_{1\leq i\leq n}\left| \left( x,y_{i}\right) \right|
\left( \sum\limits_{i,j=1}^{n}\left| \left( y_{i},y_{j}\right) \right|
\right) ^{\frac{1}{2}},\ \ x\in H.  \label{4.1}
\end{equation}
If the same choice for $c_{i}$ is made in the second inequality in (\ref{3.1}%
), then one would get 
\begin{equation*}
\left( \sum\limits_{i=1}^{n}\left| \left( x,y_{i}\right) \right| ^{2}\right)
^{2}\leq \left\| x\right\| ^{2}\max\limits_{1\leq i\leq n}\left| \left(
x,y_{i}\right) \right| \left( \sum\limits_{i=1}^{n}\left| \left(
x,y_{i}\right) \right| ^{r}\right) ^{\frac{1}{r}}\left[ \sum%
\limits_{i=1}^{n}\left( \sum\limits_{j=1}^{n}\left| \left(
y_{i},y_{j}\right) \right| \right) ^{s}\right] ^{\frac{1}{s}},
\end{equation*}
implying 
\begin{multline}
\sum\limits_{i=1}^{n}\left| \left( x,y_{i}\right) \right| ^{2}  \label{4.2}
\\
\leq \left\| x\right\| \max\limits_{1\leq i\leq n}\left| \left(
x,y_{i}\right) \right| ^{\frac{1}{2}}\left( \sum\limits_{i=1}^{n}\left|
\left( x,y_{i}\right) \right| ^{r}\right) ^{\frac{1}{2r}}\left[
\sum\limits_{i=1}^{n}\left( \sum\limits_{j=1}^{n}\left| \left(
y_{i},y_{j}\right) \right| \right) ^{s}\right] ^{\frac{1}{2s}},
\end{multline}
where $\frac{1}{r}+\frac{1}{s}=1,$ $s>1.$

The other inequalities in (\ref{3.1}) will produce the following results,
respectively 
\begin{equation}
\sum\limits_{i=1}^{n}\left| \left( x,y_{i}\right) \right| ^{2}\leq \left\|
x\right\| \max\limits_{1\leq i\leq n}\left| \left( x,y_{i}\right) \right| ^{%
\frac{1}{2}}\left( \sum\limits_{i=1}^{n}\left| \left( x,y_{i}\right) \right|
\right) ^{\frac{1}{2}}\left[ \max\limits_{1\leq i\leq n}\left(
\sum\limits_{j=1}^{n}\left| \left( y_{i},y_{j}\right) \right| \right) \right]
;  \label{4.3}
\end{equation}
\begin{multline}
\sum\limits_{i=1}^{n}\left| \left( x,y_{i}\right) \right| ^{2}  \label{4.4}
\\
\leq \left\| x\right\| \max\limits_{1\leq i\leq n}\left| \left(
x,y_{i}\right) \right| ^{\frac{1}{2}}\left( \sum\limits_{i=1}^{n}\left|
\left( x,y_{i}\right) \right| ^{p}\right) ^{\frac{1}{2p}}\left[
\sum\limits_{i=1}^{n}\left( \sum\limits_{j=1}^{n}\left| \left(
y_{i},y_{j}\right) \right| ^{q}\right) ^{\frac{1}{q}}\right] ^{\frac{1}{2}},
\end{multline}
where $p>1,$ $\frac{1}{p}+\frac{1}{q}=1;$%
\begin{multline}
\sum\limits_{i=1}^{n}\left| \left( x,y_{i}\right) \right| ^{2}  \label{4.5}
\\
\leq \left\| x\right\| \left( \sum\limits_{i=1}^{n}\left| \left(
x,y_{i}\right) \right| ^{p}\right) ^{\frac{1}{2p}}\left(
\sum\limits_{i=1}^{n}\left| \left( x,y_{i}\right) \right| ^{t}\right) ^{%
\frac{1}{2t}}\left[ \sum\limits_{i=1}^{n}\left( \sum\limits_{j=1}^{n}\left|
\left( y_{i},y_{j}\right) \right| ^{q}\right) ^{\frac{u}{q}}\right] ^{\frac{1%
}{2u}},
\end{multline}
where $p>1,$ $\frac{1}{p}+\frac{1}{q}=1,$ $t>1,$ $\frac{1}{t}+\frac{1}{u}=1;$%
\begin{multline}
\sum\limits_{i=1}^{n}\left| \left( x,y_{i}\right) \right| ^{2}  \label{4.6}
\\
\leq \left\| x\right\| \left( \sum\limits_{i=1}^{n}\left| \left(
x,y_{i}\right) \right| ^{p}\right) ^{\frac{1}{2p}}\left(
\sum\limits_{i=1}^{n}\left| \left( x,y_{i}\right) \right| \right) ^{\frac{1}{%
2}}\max\limits_{1\leq i\leq n}\left\{ \left( \sum\limits_{j=1}^{n}\left|
\left( y_{i},y_{j}\right) \right| ^{q}\right) ^{\frac{1}{2q}}\right\} ,
\end{multline}
where $p>1,$ $\frac{1}{p}+\frac{1}{q}=1;$%
\begin{equation}
\sum\limits_{i=1}^{n}\left| \left( x,y_{i}\right) \right| ^{2}\leq \left\|
x\right\| \left[ \sum\limits_{i=1}^{n}\left| \left( x,y_{i}\right) \right| %
\right] ^{\frac{1}{2}}\max\limits_{1\leq i\leq n}\left| \left(
x,y_{i}\right) \right| ^{\frac{1}{2}}\left( \sum\limits_{i=1}^{n}\left[
\max\limits_{1\leq j\leq n}\left| \left( y_{i},y_{j}\right) \right| \right]
\right) ^{\frac{1}{2}};  \label{4.7}
\end{equation}
\begin{equation}
\sum\limits_{i=1}^{n}\left| \left( x,y_{i}\right) \right| ^{2}\leq \left\|
x\right\| \left[ \sum\limits_{i=1}^{n}\left| \left( x,y_{i}\right) \right|
^{m}\right] ^{\frac{1}{2m}}\left[ \sum\limits_{i=1}^{n}\left[
\max\limits_{1\leq j\leq n}\left| \left( y_{i},y_{j}\right) \right| ^{l}%
\right] \right] ^{\frac{1}{2l}},  \label{4.8}
\end{equation}
where $m>1,$ $\frac{1}{m}+\frac{1}{l}=1;$ and 
\begin{equation}
\sum\limits_{i=1}^{n}\left| \left( x,y_{i}\right) \right| ^{2}\leq \left\|
x\right\| \sum\limits_{i=1}^{n}\left| \left( x,y_{i}\right) \right|
\max\limits_{i,1\leq j\leq n}\left| \left( y_{i},y_{j}\right) \right| ^{%
\frac{1}{2}}.  \label{4.9}
\end{equation}

If in the above inequalities we assume that $\left( y_{i}\right) _{1\leq
i\leq n}=\left( e_{i}\right) _{1\leq i\leq n},$ where $\left( e_{i}\right)
_{1\leq i\leq n}$ are orthonormal vectors in the inner product space $\left(
H,\left( \cdot ,\cdot \right) \right) ,$ then from (\ref{4.1}) -- (\ref{4.9}%
) we may deduce the following inequalities similar in a sense with Bessel's
inequality: 
\begin{equation}
\sum\limits_{i=1}^{n}\left| \left( x,e_{i}\right) \right| ^{2}\leq \sqrt{n}%
\left\| x\right\| \max\limits_{1\leq i\leq n}\left\{ \left| \left(
x,e_{i}\right) \right| \right\} ;  \label{4.10}
\end{equation}
\begin{equation}
\sum\limits_{i=1}^{n}\left| \left( x,e_{i}\right) \right| ^{2}\leq n^{\frac{1%
}{2s}}\left\| x\right\| \max\limits_{1\leq i\leq n}\left\{ \left| \left(
x,e_{i}\right) \right| ^{\frac{1}{2}}\right\} \left(
\sum\limits_{i=1}^{n}\left| \left( x,e_{i}\right) \right| ^{r}\right) ^{%
\frac{1}{2r}},  \label{4.11}
\end{equation}
where $r>1,$ $\frac{1}{r}+\frac{1}{s}=1;$%
\begin{equation}
\sum\limits_{i=1}^{n}\left| \left( x,e_{i}\right) \right| ^{2}\leq \left\|
x\right\| \max\limits_{1\leq i\leq n}\left\{ \left| \left( x,e_{i}\right)
\right| ^{\frac{1}{2}}\right\} \left( \sum\limits_{i=1}^{n}\left| \left(
x,e_{i}\right) \right| \right) ^{\frac{1}{2}};  \label{4.12}
\end{equation}
\begin{equation}
\sum\limits_{i=1}^{n}\left| \left( x,e_{i}\right) \right| ^{2}\leq \sqrt{n}%
\left\| x\right\| \max\limits_{1\leq i\leq n}\left\{ \left| \left(
x,e_{i}\right) \right| ^{\frac{1}{2}}\right\} \left(
\sum\limits_{i=1}^{n}\left| \left( x,e_{i}\right) \right| ^{p}\right) ^{%
\frac{1}{2p}},  \label{4.13}
\end{equation}
where $p>1;$%
\begin{equation}
\sum\limits_{i=1}^{n}\left| \left( x,e_{i}\right) \right| ^{2}\leq n^{\frac{1%
}{2u}}\left\| x\right\| \left( \sum\limits_{i=1}^{n}\left| \left(
x,e_{i}\right) \right| ^{p}\right) ^{\frac{1}{2p}}\left(
\sum\limits_{i=1}^{n}\left| \left( x,e_{i}\right) \right| ^{t}\right) ^{%
\frac{1}{2t}},  \label{4.14}
\end{equation}
where $p>1,$ $t>1,$ $\frac{1}{t}+\frac{1}{u}=1;$%
\begin{equation}
\sum\limits_{i=1}^{n}\left| \left( x,e_{i}\right) \right| ^{2}\leq \left\|
x\right\| \left( \sum\limits_{i=1}^{n}\left| \left( x,e_{i}\right) \right|
^{p}\right) ^{\frac{1}{2p}}\left( \sum\limits_{i=1}^{n}\left| \left(
x,e_{i}\right) \right| \right) ^{\frac{1}{2}},\ \ p>1;  \label{4.15}
\end{equation}
\begin{equation}
\sum\limits_{i=1}^{n}\left| \left( x,e_{i}\right) \right| ^{2}\leq \sqrt{n}%
\left\| x\right\| \left( \sum\limits_{i=1}^{n}\left| \left( x,e_{i}\right)
\right| \right) ^{\frac{1}{2}}\max\limits_{1\leq i\leq n}\left\{ \left|
\left( x,e_{i}\right) \right| ^{\frac{1}{2}}\right\} ;  \label{4.16}
\end{equation}
\begin{equation}
\sum\limits_{i=1}^{n}\left| \left( x,e_{i}\right) \right| ^{2}\leq n^{\frac{1%
}{2l}}\left\| x\right\| \left[ \sum\limits_{i=1}^{n}\left| \left(
x,e_{i}\right) \right| ^{m}\right] ^{\frac{1}{m}},\ \ \ m>1,\ \frac{1}{m}+%
\frac{1}{l}=1;  \label{4.17}
\end{equation}
\begin{equation}
\sum\limits_{i=1}^{n}\left| \left( x,e_{i}\right) \right| ^{2}\leq \left\|
x\right\| \sum\limits_{i=1}^{n}\left| \left( x,e_{i}\right) \right| .
\label{4.18}
\end{equation}

Corollaries \ref{c3.2} -- \ref{c3.5} will produce the following results
which do not contain the Fourier coefficients in the right side of the
inequality.

Indeed, if one chooses $c_{i}=\overline{\left( x,y_{i}\right) }$ in (\ref
{3.3}), then 
\begin{equation*}
\left( \sum\limits_{i=1}^{n}\left\vert \left( x,y_{i}\right) \right\vert
^{2}\right) ^{2}\leq \left\Vert x\right\Vert ^{2}n^{\frac{1}{p}+\frac{1}{t}%
-1}\sum\limits_{i=1}^{n}\left\vert \left( x,y_{i}\right) \right\vert ^{2}%
\left[ \sum\limits_{i=1}^{n}\left( \sum\limits_{j=1}^{n}\left\vert \left(
y_{i},y_{j}\right) \right\vert ^{q}\right) ^{\frac{u}{q}}\right] ^{\frac{1}{u%
}},
\end{equation*}
giving the following Bombieri type inequality: 
\begin{equation}
\sum\limits_{i=1}^{n}\left\vert \left( x,y_{i}\right) \right\vert ^{2}\leq
n^{\frac{1}{p}+\frac{1}{t}-1}\left\Vert x\right\Vert ^{2}\left[
\sum\limits_{i=1}^{n}\left( \sum\limits_{j=1}^{n}\left\vert \left(
y_{i},y_{j}\right) \right\vert ^{q}\right) ^{\frac{u}{q}}\right] ^{\frac{1}{u%
}},  \label{4.19}
\end{equation}
where $1<p\leq 2,$ $1<t\leq 2,$ $\frac{1}{p}+\frac{1}{q}=1,$ $\frac{1}{t}+%
\frac{1}{u}=1.$

If in this inequality we consider $p=q=t=u=2,$ then 
\begin{equation}
\sum\limits_{i=1}^{n}\left| \left( x,y_{i}\right) \right| ^{2}\leq \left\|
x\right\| ^{2}\left( \sum\limits_{i,j=1}^{n}\left| \left( y_{i},y_{j}\right)
\right| ^{2}\right) ^{\frac{1}{2}}.  \label{4.20}
\end{equation}
For a different proof of this result see also \cite{5b}.

In a similar way, if $c_{i}=\overline{\left( x,y_{i}\right) }$ in (\ref{3.6}%
), then 
\begin{equation}
\sum\limits_{i=1}^{n}\left| \left( x,y_{i}\right) \right| ^{2}\leq n^{\frac{1%
}{m}}\left\| x\right\| ^{2}\left( \sum\limits_{i=1}^{n}\left[
\max\limits_{1\leq j\leq n}\left| \left( y_{i},y_{j}\right) \right| \right]
^{l}\right) ^{\frac{1}{l}},  \label{4.21}
\end{equation}
where $m>1,$ $\frac{1}{m}+\frac{1}{l}=1.$

Finally, if $c_{i}=\overline{\left( x,y_{i}\right) }$ $\left( i=1,\dots
,n\right) ,$ is taken in (\ref{3.7}), then 
\begin{equation}
\sum\limits_{i=1}^{n}\left| \left( x,y_{i}\right) \right| ^{2}\leq n\left\|
x\right\| ^{2}\max\limits_{1\leq i,j\leq n}\left| \left( y_{i},y_{j}\right)
\right| .  \label{4.22}
\end{equation}

\begin{remark}
\label{r4.1}Let us compare Bombieri's result 
\begin{equation}
\sum\limits_{i=1}^{n}\left| \left( x,y_{i}\right) \right| ^{2}\leq \left\|
x\right\| ^{2}\max\limits_{1\leq i\leq n}\left\{ \sum\limits_{j=1}^{n}\left|
\left( y_{i},y_{j}\right) \right| \right\}  \label{4.23}
\end{equation}
with our result 
\begin{equation}
\sum\limits_{i=1}^{n}\left| \left( x,y_{i}\right) \right| ^{2}\leq \left\|
x\right\| ^{2}\left\{ \sum\limits_{i,j=1}^{n}\left| \left(
y_{i},y_{j}\right) \right| ^{2}\right\} ^{\frac{1}{2}}.  \label{4.24}
\end{equation}
Denote 
\begin{equation*}
M_{1}:=\max_{1\leq i\leq n}\left\{ \sum_{j=1}^{n}\left| \left(
y_{i},y_{j}\right) \right| \right\}
\end{equation*}
and 
\begin{equation*}
M_{2}:=\left[ \sum_{i,j=1}^{n}\left| \left( y_{i},y_{j}\right) \right| ^{2}%
\right] ^{\frac{1}{2}}.
\end{equation*}

If we choose the inner product space $H=\mathbb{R}$, $\left( x,y\right) :=xy$
and $n=2,$ then for $y_{1}=a,$ $y_{2}=b,$ $a,b>0,$ we have 
\begin{gather*}
M_{1}=\max \left\{ a^{2}+ab,ab+b^{2}\right\} =\left( a+b\right) \max \left(
a,b\right) , \\
M_{2}=\left( a^{4}+a^{2}b^{2}+a^{2}b^{2}+b^{4}\right) ^{\frac{1}{2}%
}=a^{2}+b^{2}.
\end{gather*}

Assume that $a\geq b.$ Then $M_{1}=a^{2}+ab\geq a^{2}+b^{2}=M_{2},$ showing
that, in this case, the bound provided by (\ref{4.24}) is better than the
bound provided by (\ref{4.23}). If $\left( y_{i}\right) _{1\leq i\leq n}$
are orthonormal vectors, then $M_{1}=1,$ $M_{2}=\sqrt{n}$, showing that in
this case the Bombieri inequality (which becomes Bessel's inequality)
provides a better bound than (\ref{4.24}).
\end{remark}


\begin{thebibliography}{9}
\bibitem{1b}  R. BELLMAN, Almost orthogonal series, \textit{Bull. Amer.
Math. Soc., }\textbf{50} (1944), 517--519.

\bibitem{2b}  R.P. BOAS, A general moment problem, \textit{Amer. J. Math., }%
\textbf{63} (1941), 361--370.

\bibitem{2ab}  E. BOMBIERI, A note on the large sieve, \textit{Acta Arith.}, 
\textbf{18}(1971), 401-404.

\bibitem{3b}  S.S. DRAGOMIR and J. S\'{A}NDOR, On Bessel's and Gaur's
inequality in prehilbertian spaces, \textit{Periodica Math. Hung., }\textbf{%
29}(3) (1994), 197--205.

\bibitem{4b}  S.S. DRAGOMIR and B. MOND, On the Boas-Bellman generalisation
of Bessel's inequality in inner product spaces, \textit{Italian J. of Pure
\& Appl. Math., }\textbf{3} (1998), 29--35.

\bibitem{5b}  S.S. DRAGOMIR, B. MOND and J.E. PE\v{C}ARI\'{C}, Some remarks
on Bessel's inequality in inner product spaces, \textit{Studia Univ. Babe\c{s%
}-Bolyai, Mathematica, }\textbf{37}(4) (1992), 77--86.

\bibitem{5ab}  H. HEILBRONN, On the averages of some arithmetical functions
of two variables, \textit{Mathematica}, \textbf{5}(1958), 1-7.

\bibitem{6b}  D.S. MITRINOVI\'{C}, J.E. PE\v{C}ARI\'{C} and A.M. FINK, 
\textit{Classical and New Inequalities in Analysis, }Kluwer Academic
Publishers, 1993.

\bibitem{7b}  J.E. PE\v{C}ARI\'{C}, On some classical inequalities in
unitary spaces, \textit{Mat. Bilten} (Scopje), \textbf{16}(1992), 63-72.
\end{thebibliography}
\end{document}